\documentclass[12pt]{article}
\usepackage{latexsym,amstex}
\usepackage[all]{xy}
\xyoption{v2}
\define{\scA}{{\cal A}}
\define{\scB}{{\cal B}}
\define{\scC}{{\cal C}}
\define{\scD}{{\cal D}}
\define{\scK}{{\cal K}}
\define{\scL}{{\cal L}}
\define{\scM}{{\cal M}}
\define{\scP}{{\cal P}}
\define{\gtensor}{{\hat{\otimes}}}

\begin{document}
\centerline{\large\bf Constructions of $E_n$ Operads}
\centerline{Z. Fiedorowicz}

\bigskip
Throughout this talk I will use the following conventions and notations.
I will primarily consider operads in the category of compactly generated
Hausdorff topological spaces having the homotopy type of $CW$-complexes.
When I refer to simplicial operads or operads in the category of posets,
it will be understood that they can be converted to topological operads
by taking geometric realization or geometric realization of the nerve,
respectively. I will only consider operads $\scA$ satisfying $\scA(0)=*$.
This allows us to define degeneracy maps
$$d_i:\scA(n)\longrightarrow\scA(n-1)\quad i=1,2,\dots,n$$
by composing with $*\in\scA(0)$. We will refer to collections of
spaces $\scA=\{\scA(k)\}_{k\ge 0}$ equipped with symmetric group actions
and such degeneracy maps, satisfying the evident relations, but without
any further operad composition maps, as {\it preoperads\/}. (We refer the
reader to Berger's talk in this volume \cite{Be2} (or \cite{Be}) for a
precise definition.)

\bigskip
\noindent
{\bf Definition.} An operad map $\scA\to\scB$ is said to be an
{\it equivalence\/} if for each $k\ge 0$, $\scA(k)\to\scB(k)$ is a
$\Sigma_k$-equivariant homotopy equivalence. We say that an operad $\scA$ is
$E_n$ ($n=1,2,3,\dots,\infty$) if there is a chain of operad equivalences
(in either or both directions) connecting $\scA$ to the Boardman-Vogt little
$n$-cubes operad $\scC_n$ (cf. \cite{BV}).

\bigskip
\noindent
{\bf Remark.} It can be shown that an operad $\scA$ is $E_n$ iff there
is an operad equivalence $W\scC_n\to\scA$, where $W\scC_n$ is a certain
enlarged cofibrant model of the little $n$-cubes operad. (See Vogt's talk
in this volume \cite{RV} (or \cite{BV}) for details.) However this is of
little use in recognizing $E_n$ operads, since the only practical way to
construct such an equivalence $W\scC_n\to\scA$ is to be given a chain of
operad equivalences connecting $\scA$ to $\scC_n$.

\vskip0.3in

The following results provide useful criteria for recognizing $E_n$ operads.

\bigskip
\noindent
{\bf Theorem (Recognition principles for $E_n$ operads.)}
\begin{enumerate}
\item $\scA$ is an $E_\infty$ operad iff
\begin{description}
\item[{\rm (i)}] Each space $\scA(k)$ is contractible.
\item[{\rm (ii)}] $\Sigma_k$ acts freely on $\scA(k)$ for each $k\ge 0$.
\end{description}
\item $\scA$ is an $E_1$ (a.k.a. $A_\infty$) operad iff
\begin{description}
\item[{\rm (i)}] Each path component of each space $\scA(k)$ is contractible.
\item[{\rm (ii)}] $\Sigma_k$ acts freely and transitively on $\pi_0\scA(k)$
for each $k\ge 0$.
\end{description}
\item $\scA$ is an $E_2$ operad iff each space $\scA(k)$ is connected and
the collection of universal covering spaces $\{\widetilde{\scA}(k)\}_{k\ge 0}$
forms a $B_\infty$ operad. That is
\begin{description}
\item[{\rm (i)}] Each space $\widetilde{\scA}(k)$ is contractible.
\item[{\rm (ii)}] The braid group $B_k$ acts freely on $\widetilde{\scA}(k)$
for each $k\ge 0$.
\end{description}
\end{enumerate}

\bigskip
{\it Proof Sketch.\/} The first two results are classical, due to
Boardman-Vogt \cite{BV}.
The fact that $E_\infty$, resp. $E_1$, operads satisfy these criteria is
immediate from the fact that $\scC_\infty$, resp. $\scC_1$, satisfy these
criteria, and that these criteria are preserved by operad equivalences.
For the reverse implications, one notes that the given criteria insure
that the projections provide a chain of operad equivalences
$$\scA\stackrel{\simeq}{\longleftarrow}\scA\times\scC_\infty
\stackrel{\simeq}{\longrightarrow}\scC_\infty$$
in the $E_\infty$ case, and a similar chain of operad equivalences
$$\scA\stackrel{\simeq}{\longleftarrow}\scA\times_{\pi_0(\scA)}\scC_1
\stackrel{\simeq}{\longrightarrow}\scC_1$$
in the $E_1$ case.

The recognition principle for $E_2$ operads is due to the author (cf.
\cite{ZF}). It first requires noting that the definition of operad can be
reformulated using actions by braid groups in place of actions by symmetric
groups. One then notes that the operad structure of an $E_2$ operad $\scA$ can
be lifted to a (braided) operad structure on the universal covers
$\{\widetilde{\scA}\}_{k\ge 0}$. The lifting requires a consistent
choice of basepoints in $\scA$. Such a choice is provided by mapping
an $E_1$ operad into $\scA$ as follows:
$$W\scC_1\longrightarrow W\scC_2\stackrel{\simeq}{\longrightarrow}\scA,$$
(cf. remark above). This gives the requisite $B_\infty$ structure on
$\widetilde{\scA}$. Conversely given a braided operad structure on
$\widetilde{\scA}$, one can recover a regular operad structure on $\scA$,
by taking the orbit spaces of the actions by the pure braid groups $PB_k$. The
rest of the proof now proceeds as in the $E_\infty$ case:
$$\scA\stackrel{\simeq}{\longleftarrow}\left.\left(\widetilde{\scA}\times
\widetilde{\scC}_2\right)\right/PB_*\stackrel{\simeq}{\longrightarrow}\scC_2.$$

\bigskip
Unfortunately at this time, there is nothing like a recognition principle
for $E_n$ operads for $3\le n<\infty$.  Indeed until the last five years
or so, the only known $E_n$ operads in such cases were minor variants
of the little $n$-cubes operads, eg. little $n$-disks operads.

This situation changed with the discovery of a nice family of simplicial
$E_n$ operads and operad equivalences (for all $n$):

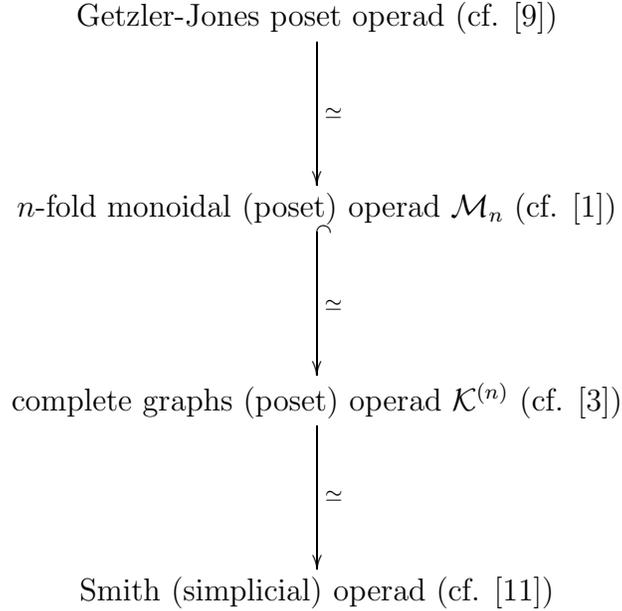
$$\xymatrix{
\mbox{Getzler-Jones poset operad (cf. \cite{GJ})}\ar[dd]^{\simeq}\\ \\
\mbox{$n$-fold monoidal (poset) operad\ }\scM_n\mbox{ (cf. \cite{BFSV})}
\ar@{^{(}->}[dd]^{\simeq}\\ \\
\mbox{complete graphs (poset) operad\ }\scK^{(n)}\mbox{ (cf. \cite{Be})}
\ar[dd]^{\simeq}\\ \\
\mbox{Smith (simplicial) operad (cf. \cite{Sm})}
}$$

I will not discuss these examples and instead refer the reader to
Berger's talk in this volume \cite{Be2} for more details. However I would
like to briefly discuss the method developed by Berger to prove that these
examples are $E_n$ operads.

This method is the notion of a cellular decomposition of a topological operad
over a poset operad. Let me first discuss the simpler notion of a cellular
decomposition of a topological space $X$ over a finite poset $\scP$. By
this we mean a decomposition $X=\cup_{P\in\scP}C_P$, where each cell $C_P$ is
a closed subset of $X$. We further suppose that $P<Q\Longrightarrow C_P
\subseteq C_Q$ and define $\mbox{bd}(C_P)=\cup_{P'<P}C_{P'}$,
$\mbox{int}(C_P)=C_P-\mbox{bd}(C_P)$. If we assume in addition that
\begin{enumerate}
\item Each $C_P$ is contractible.
\item The inclusions $\mbox{bd}(C_P)\subseteq C_P$ are cofibrations.
\item $\mbox{int}(C_P)\cap\mbox{int}(C_Q)=\emptyset$ if $P\ne Q$.
\end{enumerate}
Then we have a sequence of equivalences and homeomorphisms:
$$|\scP|\stackrel{\simeq}{\longleftarrow}\mbox{hocolim}_{P\in\scP}C_P
\stackrel{\simeq}{\longrightarrow}\mbox{colim}_{P\in\scP}C_P
\stackrel{\cong}{\longrightarrow}X,$$
where $|\scP|$ denotes the geometric realization of the nerve of $\scP$.
Briefly, we need the hypothesis (1) to insure that the first arrow is
an equivalence, hypothesis (2) to insure that the second arrow is an
equivalence, and hypothesis (3) to insure that the third arrow is a
homeomorphism. It is useful to note that if $X$ has a cellular decomposition
over a poset $\scP$, then the product $X^m$ has a cellular decomposition
indexed over the product poset $\scP^m$.

Please note that the word ``cellular'' in the above definition is being
used in a very loose sense. It is {\bf not} to be supposed that either
the cells or their interiors look anything like disks. We merely assume
that the cells satisfy the above three conditions.

To upgrade the notion of cellular decomposition to the context of operads,
we assume that the posets $\{\scP(k)\}_{k\ge 0}$ indexing the cellular
decompositions of the topological operad $\scA$, themselves form an operad
in the category of posets.  In addition we assume the evident compatibility
conditions between the operad structures of $\scP$ and $\scA$:
\begin{enumerate}
\item $1_\scA\in C_{1_\scP}$, where $1_\scA\in\scA(1)$, $1_\scP\in\scP(1)$
are the units of the operads $\scA$, $\scP$ respectively.
\item $(C_P)\sigma=C_{P\sigma}$, for any $\sigma\in\Sigma_k$, $P\in\scP(k)$.
\item If $\mu_\scA$, $\mu_\scP$ denote the operad composition maps in
$\scA$, $\scP$ respectively, then
$$\mu_\scA\left(C_{P}\times C_{Q_1}\times C_{Q_2}\times\dots C_{Q_k}\right)\subseteq C_{\mu_\scP(P;Q_1,Q_2,\dots,Q_k)}$$
\end{enumerate}
Under these hypotheses the sequence:
$$\{|\scP(k)|\}_{k\ge 0}\stackrel{\simeq}{\longleftarrow}
\{\mbox{hocolim}_{P\in\scP}C_P(k)\}_{k\ge 0}
\stackrel{\simeq}{\longrightarrow}\{\mbox{colim}_{P\in\scP(k)}C_P\}_{k\ge 0}
\stackrel{\cong}{\longrightarrow}\{\scA(k)\}_{k\ge 0}$$
is a chain of operad equivalences. Thus under these circumstances
the topological operad $\scA$ is $E_n$ iff the poset operad $\scP$ is $E_n$.
Currently this is the only practical method for recognizing $E_n$ operads
for $3\le n<\infty$.

\bigskip
\noindent
{\bf Remark.} The Getzler-Jones poset operad has not been explicitly
described in the literature to our knowledge.  However Getzler-Jones
describe a cellular decomposition of their topological operad given
by compactifications of configuration spaces, and it is not difficult
to see that their decomposition is indexed by a poset operad which maps
into the $n$-fold monoidal operad, and it can be shown that this map of
poset operads is an equivalence. This also furnishes a proof that the
Getzler-Jones topological operad is $E_n$. (There may be unpublished
more direct proofs of this fact, but we are unaware of the details.)

\bigskip
I recently noticed that this method can be combined with the following
right adjoint construction to produce many new examples of $E_n$ operads.

\bigskip
\noindent
{\bf Theorem 1.} (i) The forgetful functor
$$\xymatrix{
U:\ \mbox{Preoperads}\ar[r] &\Bbb Z/2\mbox{\ -spaces}\\
\scA \ar@{|->}[r] &\scA(2)
}$$
has a right adjoint
$$R:\ \Bbb Z/2\mbox{\ -spaces}\longrightarrow \mbox{Preoperads}$$\newline
(ii) For any $\Bbb Z/2$-space $X$, $RX$ has a functorial operad structure.
\newline
(iii) If $\scA$ is an operad such that $\scA(1)=\{1_\scA\}$, then
$$\scA\longrightarrow RU\scA=R\left(\scA(2)\right)$$
is a map of operads.

\bigskip
\noindent
{\bf Remark.} Without the hypothesis that $\scA(1)=\{1_\scA\}$
$$\scA\longrightarrow RU\scA=R\left(\scA(2)\right)$$
is only a map of preoperads, ie. it is only compatible with the symmetric
group actions and composition with constants, {\bf not} with general
compositions. We shall illustrate this phenomenon below with the little
$n$-cubes operad.  I should add that I initially overlooked the necessity
of this condition and noticed it only while preparing this talk.

\bigskip
{\it Proof Sketch.\/} Given a $\Bbb Z/2$-space $X$, we need to construct
a sequence of spaces $\{RX(k)\}_{k\ge 0}$ with an appropriate structure.
As a topological space $RX(k)$ is just the product space $X^{\binom{k}{2}}$,
where $\binom{k}{2}$ is the binomial coefficient $k(k-1)/2$. However to
describe the symmetric group actions and the operad compositions, it is
convenient to describe $RX(k)$ as follows.

Consider the complete graph on the set of vertices $\{1,2,\dots,k\}$.
By an orientation of this graph we mean the assignment of a direction to
each edge of the graph.  By a labelling of this graph, we mean the
assignment of a point in $X$ to each edge of the graph. (Different edges
may labelled by different points.) We further impose the equivalence
relation that changing the direction of an edge is equivalent to changing
its label from $x$ to $\overline{x}$, where $x\mapsto\overline{x}$ denotes
the $\Bbb Z/2$ action on $X$. We define $RX(k)$ as the space of all
such orientations and labellings of the complete graph on $\{1,2,\dots,k\}$.
Formally we can describe the elements of $RX(k)$ as functions $f(i,j)$,
where $i\ne j\ \in\ \{1,2,\dots,k\}$, taking values in $X$, and satisfying
the condition that $f(j,i)=\overline{f(i,j)}$.  $RX(k)$ can be identified
with the product $X^{\binom{k}{2}}$ by noting that using the equivalence
relation, we can choose  unique representatives in $RX(k)$ with the edges
oriented in the direction of the natural order on $\{1,2,\dots,k\}$, and with
an arbitrary labelling of the edges.

The symmetric group action on $RX(k)$ is induced by the natural action of
$\Sigma_k$ on the vertex set $\{1,2,\dots,k\}$. The operad composition
$$RX(k)\times RX(i_1)\times RX(i_2)\times\dots\times RX(i_k)
\longrightarrow RX(i_1+i_2+\dots+i_k)$$
is also easily described with the complete graph formalism. We need to
describe an orientation and labelling of the complete graph on the set
$\{1,2,\dots,i_1+i_2+\dots+i_k\}$ using the orientation and labelling
information encoded in the product space on the left. It is almost
selfevident how to do this: for edges connecting vertices in a given
block 
$$\{(i_1+i_2+\dots+i_{j-1})+1, (i_1+i_2+\dots+i_{j-1})+2,\dots,
i_1+i_2+\dots+i_{j-1}+i_j\}$$
use the orientation and labelling information from the space $RX(i_j)$.
For edges connecting vertices in different blocks use the orientation
and labelling information encoded in the space $RX(k)$. The fact that
our construction specifies an operad is then an exercise in understanding
the notation.

For a preoperad $\scA$ the unit of the adjunction
$$\scA(k)\longrightarrow (RU\scA)(k)=\left(\scA(2)\right)^{\binom{k}{2}}$$
is specified by the obvious iterated degeneracies. It is obvious that
$RU$ specifies an idempotent monad on the category of preoperads
(since $URU=U$). Hence it follows that $R$ is right adjoint to $U$.

The proof of part (iii) is deferred for now. It will be proved below
using a more structured right adjoint construction $R_2$ which coincides
with $R$ under the hypothesis that $\scA(1)=\{1_\scA\}$.

\bigskip
\noindent
{\bf Example.} The little $2$-cubes operad $\scC_2$ does not satisfy the
hypothesis of part (iii) of the theorem: $\scC_2(1)\ne\{1_{\scC_2}\}$.
We would like to show that
$$\scC_2\longrightarrow RU\scC_2$$
is not an operad map by showing that the diagram below does {\bf not} commute
$$\xymatrix{
\scC_2(2)\times\scC_2(2)\times\scC_2(1)\ar[r]\ar[d]
&\scC_2(3)\ar[d]\\
(RU\scC_2)(2)\times(RU\scC_2)(2)\times(RU\scC_2)(1)\ar[r]
&(RU\scC_2)(3)
}$$
where the vertical maps are given by the unit $\scC_2\to RU\scC_2$
and the horizontal maps are compositions in the operads $\scC_2$, $RU\scC_2$
respectively.

To see this we take a typical element in the top left corner and chase it
around both sides of the diagram:

\vskip-10pt
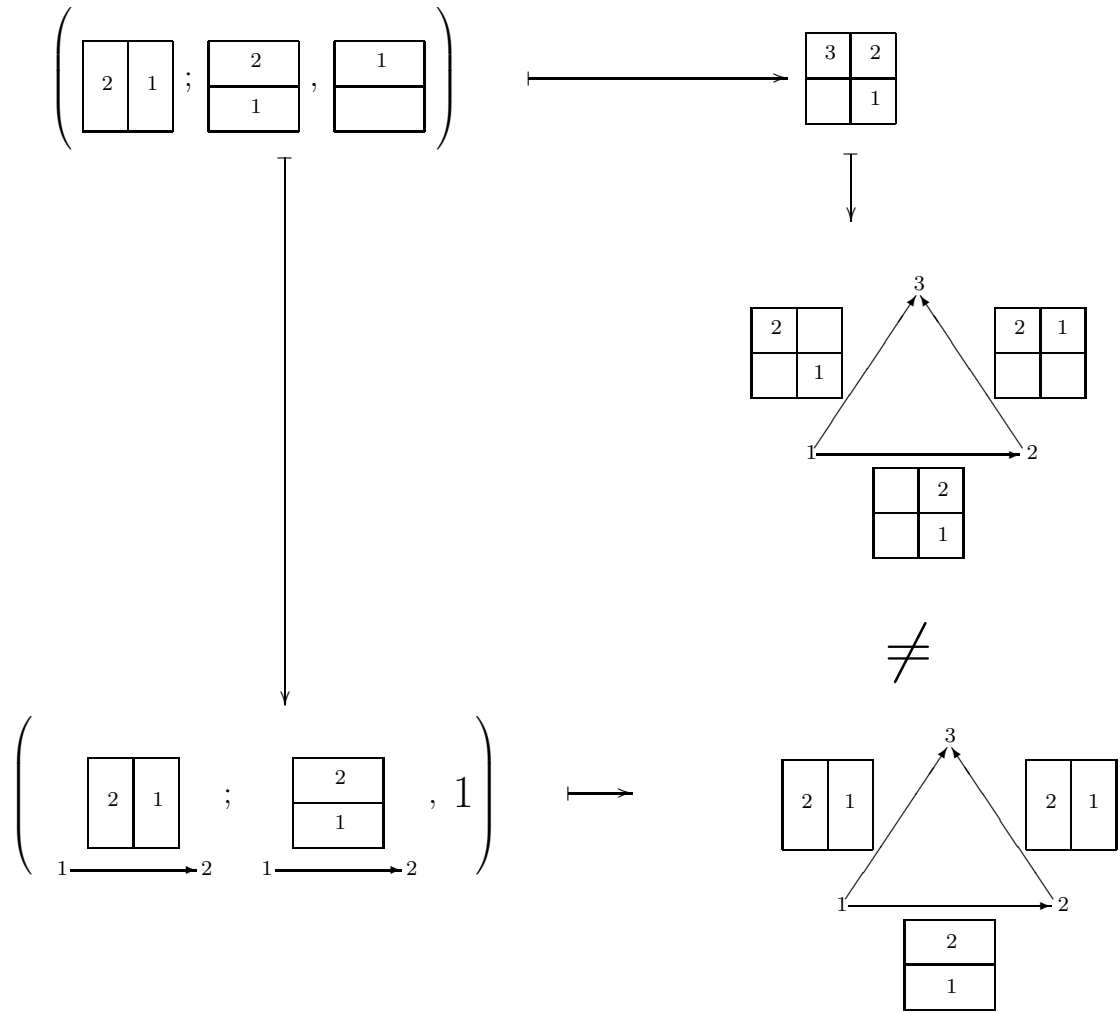
$$\xymatrix{
\left(\mbox{
\unitlength0.3cm
\begin{picture}(4.5, 3.5)
\put(0, -2) {\line(1, 0) {4}}
\put(0, -2) {\line(0, 1) {4}}
\put(0, 2) {\line(1, 0) {4}}
\put(4, -2) {\line(0, 1) {4}}
\put(2, -2) {\line(0, 1) {4}}
\put(2.85, -0.15) {$\scriptstyle 1$}
\put(0.85, -0.15) {$\scriptstyle 2$}
\end{picture}};
\mbox{
\unitlength0.3cm
\begin{picture}(4.5, 3.0)
\put(0, -2) {\line(1, 0) {4}}
\put(0, 0) {\line(1, 0) {4}}
\put(0, 2) {\line(1, 0) {4}}
\put(4, -2) {\line(0, 1) {4}}
\put(0, -2) {\line(0, 1) {4}}
\put(1.85, -1.15) {$\scriptstyle 1$}
\put(1.85, 0.85) {$\scriptstyle 2$}
\end{picture}},
\mbox{
\unitlength0.3cm
\begin{picture}(4.5, 3.0)
\put(0, -2) {\line(1, 0) {4}}
\put(0, 0) {\line(1, 0) {4}}
\put(0, 2) {\line(1, 0) {4}}
\put(4, -2) {\line(0, 1) {4}}
\put(0, -2) {\line(0, 1) {4}}
\put(1.85, 0.85) {$\scriptstyle 1$}
\end{picture}}
\right)\qquad\ar@{|->}[r]\ar@{|->}[ddd]
&\mbox{
\unitlength0.3cm
\begin{picture}(4.5, 6.0)
\put(0, 1) {\line(1, 0) {4}}
\put(0, 1) {\line(0, 1) {4}}
\put(0, 5) {\line(1, 0) {4}}
\put(4, 1) {\line(0, 1) {4}}
\put(2, 1) {\line(0, 1) {4}}
\put(0, 3) {\line(1, 0) {4}}
\put(2.85, 3.85) {$\scriptstyle 2$}
\put(0.85, 3.85) {$\scriptstyle 3$}
\put(2.85, 1.85) {$\scriptstyle 1$}
\end{picture}}\ar@{|->}[d]\\
&\qquad\qquad\mbox{
\unitlength0.3cm
\begin{picture}(10, 10)
\put(0.4, 0.3) {\vector(2, 3) {4.5}}
\put(9.6, 0.3) {\vector(-2, 3) {4.5}}
\put(0.5, 0) {\vector(1, 0) {9.0}}
\put(0, -.2) {$\scriptstyle 1$}
\put(4.8, 7.3) {$\scriptstyle 3$}
\put(9.8, -.2) {$\scriptstyle 2$}
\put(-2.4, 2.5) {\line(1,0) {4}}
\put(-2.4, 2.5) {\line(0, 1) {4}}
\put(1.6, 2.5) {\line(0, 1) {4}}
\put(-2.4, 6.5) {\line(1, 0) {4}}
\put(-0.4, 2.5) {\line(0, 1) {4}}
\put(-2.4, 4.5) {\line(1, 0) {4}}
\put(-1.55, 5.35)  {$\scriptstyle 2$}
\put(0.35, 3.35)  {$\scriptstyle 1$}
\put(8.4, 2.5) {\line(1,0) {4}}
\put(8.4, 2.5) {\line(0, 1) {4}}
\put(12.4, 2.5) {\line(0, 1) {4}}
\put(8.4, 6.5) {\line(1, 0) {4}}
\put(10.4, 2.5) {\line(0, 1) {4}}
\put(8.4, 4.5) {\line(1, 0) {4}}
\put(9.25, 5.35)  {$\scriptstyle 2$}
\put(11.15, 5.35)  {$\scriptstyle 1$}
\put(3, -0.6) {\line(1, 0) {4}}
\put(3, -0.6) {\line(0, -1) {4}}
\put(3, -4.6) {\line(1, 0) {4}}
\put(7, -0.6) {\line(0, -1) {4}}
\put(5, -0.6) {\line(0, -1) {4}}
\put(3, -2.6) {\line(1, 0) {4}}
\put(5.85, -3.85) {$\scriptstyle 1$}
\put(5.85, -1.85) {$\scriptstyle 2$}
\end{picture}}\\
&\qquad\qquad\mbox{
\unitlength0.3cm
\begin{picture}(3, 3)
\put(0,-3) {\Huge $\ne$}
\end{picture}}
\\
\left(\mbox{
\unitlength0.3cm
\begin{picture}(8, 4.0)
\put(2, -2) {\line(1, 0) {4}}
\put(2, -2) {\line(0, 1) {4}}
\put(2, 2) {\line(1, 0) {4}}
\put(6, -2) {\line(0, 1) {4}}
\put(4, -2) {\line(0, 1) {4}}
\put(4.85, -0.15) {$\scriptstyle 1$}
\put(2.85, -0.15) {$\scriptstyle 2$}
\put(1.2,-3) {\vector(1,0) {5.6}}
\put(0.6,-3.2) {$\scriptstyle 1$}
\put(7.0,-3.2) {$\scriptstyle 2$}
\end{picture}};
\mbox{
\unitlength0.3cm
\begin{picture}(8, 4.0)
\put(2, -2) {\line(1, 0) {4}}
\put(2, 0) {\line(1, 0) {4}}
\put(2, 2) {\line(1, 0) {4}}
\put(6, -2) {\line(0, 1) {4}}
\put(2, -2) {\line(0, 1) {4}}
\put(3.85, -1.15) {$\scriptstyle 1$}
\put(3.85, 0.85) {$\scriptstyle 2$}
\put(1.2,-3) {\vector(1,0) {5.6}}
\put(0.6,-3.2) {$\scriptstyle 1$}
\put(7.0,-3.2) {$\scriptstyle 2$}
\end{picture}},
\mbox{
\unitlength0.3cm
\begin{picture}(1, 3.0)
\put(0, -.2) {\Large 1}
\end{picture}}
\right)\qquad\ar@{|->}[r]
&\qquad\qquad\qquad\mbox{
\unitlength0.3cm
\begin{picture}(10, 10)
\put(0.4, 0.3) {\vector(2, 3) {4.5}}
\put(9.6, 0.3) {\vector(-2, 3) {4.5}}
\put(0.5, 0) {\vector(1, 0) {9.0}}
\put(0, -.2) {$\scriptstyle 1$}
\put(4.8, 7.3) {$\scriptstyle 3$}
\put(9.8, -.2) {$\scriptstyle 2$}
\put(-2.4, 2.5) {\line(1,0) {4}}
\put(-2.4, 2.5) {\line(0, 1) {4}}
\put(1.6, 2.5) {\line(0, 1) {4}}
\put(-2.4, 6.5) {\line(1, 0) {4}}
\put(-0.4, 2.5) {\line(0, 1) {4}}
\put(-1.55, 4.35)  {$\scriptstyle 2$}
\put(0.35, 4.35)  {$\scriptstyle 1$}
\put(8.4, 2.5) {\line(1,0) {4}}
\put(8.4, 2.5) {\line(0, 1) {4}}
\put(12.4, 2.5) {\line(0, 1) {4}}
\put(8.4, 6.5) {\line(1, 0) {4}}
\put(10.4, 2.5) {\line(0, 1) {4}}
\put(9.25, 4.35)  {$\scriptstyle 2$}
\put(11.15, 4.35)  {$\scriptstyle 1$}
\put(3, -0.6) {\line(1, 0) {4}}
\put(3, -0.6) {\line(0, -1) {4}}
\put(3, -4.6) {\line(1, 0) {4}}
\put(7, -0.6) {\line(0, -1) {4}}
\put(3, -2.6) {\line(1, 0) {4}}
\put(4.85, -3.85) {$\scriptstyle 1$}
\put(4.85, -1.85) {$\scriptstyle 2$}
\end{picture}}
\\ \\
}$$

\bigskip
\noindent
Thus we see that the image of $\scC_2$ in $RU\scC_2$ is not closed under
operad composition in $RU\scC_2$.

\vskip15pt
\noindent
{\bf Remark.} Theorem 1 has a poset version: the forgetful functor
$$U:\mbox{Poset preoperads}\longrightarrow\Bbb Z/2\mbox{ -posets}$$
has a right adjoint
$$R:\Bbb Z/2\mbox{ -posets}\longrightarrow\mbox{Poset preoperads}$$
satisfying the same conditions as above. The construction of the functor
$R$ is identical to the case of spaces.   If a $\Bbb Z/2$-space $X$ has a
cellular decomposition over $\Bbb Z/2$-poset $\scP$, with the $\Bbb Z/2$
actions being compatible, then the topological operad $RX$ has a cellular
decomposition over the poset operad $R\scP$. Then using Berger's argument
we have a chain of operad equivalences connecting $|R\scP|$ to $RX$.

\bigskip
For any $\Bbb Z/2$-space $X$ the operad $RX$ is not $E_n$ for at least two
reasons:
\begin{enumerate}
\item The spaces $RX(k)$ have the wrong homotopy type. For by definition
$RX(k)=X^{\binom{k}{2}}$, whereas the $k$-th space of an $E_n$ operad
can not have the homotopy type of a product of multiple copies of a space
unless $n=\infty$, $k=0$ or $k=1$.
\item For $k>2$ the action of $\Sigma_k$ on $RX(k)$ isn't free. For example
for any $x\in X$ the following point in $RX(3)$
$$\xymatrix{
&3\ar[dl]_x \\
1\ar[rr]^x &&2\ar[ul]_x
}\quad\begin{array}{c}\\ \\ \mbox{\Large $=$} \\ \end{array}\quad
\xymatrix{
&3 \\
1\ar[ur]^{\overline{x}}\ar[rr]^x &&2\ar[ul]_x
}$$
is fixed by the cyclic permutation $(1\ 2 \ 3)$.
\end{enumerate}
However under suitable conditions, $RX$ has many $E_n$ suboperads. We begin
with a definition.

\bigskip
\noindent
{\bf Definition.} Let $X$ be a space which is homotopy equivalent to $S^{n-1}$
with a free $\Bbb Z/2$ action.  If $n=1$ we say that $X$ has a
{\it hemispherical cellular decomposition\/} if the $\Bbb Z/2$ action permutes
the two path components. If $n>1$ we define the notion of a hemispherical
cellular decomposition for $X$ by recursively requiring:
\begin{description}
\item[{\rm (i)}] $X$ has a decomposition $X=D\cup\overline{D}$, where
$\overline{D}$ is the image of $D$ under the $\Bbb Z/2$ action.
\item[{\rm (ii)}] $D$ is contractible by a contraction $H:D\times I\to D$
having the property that $H(x,t)\not\in\overline{D}$ for $t>0$.
\item[{\rm (iii)}] $D\cap\overline{D}$ is homotopy equivalent to $S^{n-2}$
and has a hemispherical cellular decomposition.
\end{description}
We can think of $D\subset X$ as being a kind of fundamental domain for the
$\Bbb Z/2$ action.

It is an immediate consequence that if $X$ has a hemispherical cellular
decomposition, then $X$ has a cellular decomposition over the $\Bbb Z/2$ poset
$$\scK^{(n)}(2)=\{1,2,\dots, n\}\times\Bbb Z/2$$
with partial order
$$(i,\sigma)<(j,\tau)\quad\mbox{if }i<j$$
and with $\Bbb Z/2$ action via the second factor. As suggested by the notation,
$\scK^{(n)}(2)$ is the 2-space of Berger's complete graphs operad $\scK^{(n)}$.
Now by the remark above, the topological operad $RX$ has a cellular
decomposition over the poset operad $R\left(\scK^{(n)}(2)\right)=RU\scK^{(n)}$.
Now for every poset suboperad $\scP\subset RU\scK^{(n)}$ we define
$$\left(R_{\scP}X\right)(k)=\cup_{P\in\scP(k)}\mbox{int}C_P\subset RX(k)$$
and we denote $R_{\scP}X=
\left\{\left(R_{\scP}(k)X\right)\right\}_{k\ge 0}$.

\bigskip
Our main result is

\bigskip
\noindent
{\bf Theorem 2.} Let $X$ be a space homotopy equivalent to $S^{n-1}$,
with a free $\Bbb Z/2$ action and a hemispherical cellular decomposition.
Let $\scP$ be an $E_n$ poset suboperad of $RU\scK^{(n)}$
(eg. $\scP=\scK^{(n)}$ or $\scP=\scM_n$). Then $R_{\scP}X$ is an $E_n$
suboperad of $RX$.

\bigskip
{\it Proof Sketch.\/} It is almost clear that $R_{\scP}X$ has a cellular
decomposition over the poset operad $\scP$. The only nonobvious thing to
check is that the cells of the decomposition are contractible. For the
cells of $R_{\scP}X$ are obtained from the cells of $RX$ by removing
boundary cells indexed by elements outside the poset $\scP$. Their
contractibility is insured by condition (ii) in the definition above.
(One also has to appeal to results of \cite{Do} to check that the required
cofibration conditions for a cellular decomposition continue to hold in
$R_{\scP}X$.) Now we can apply Berger's argument to obtain a chain of
operad equivalences connecting $R_{\scP}X$ to the $E_n$ operad $|\scP|$.

\bigskip
\noindent
{\bf Remark.} The poset operad $RU\scK^{(n)}$ has many other $E_n$ suboperads
besides the $n$-fold monoidal operad $\scM_n$ and the complete graphs operad
$\scK^{(n)}$. It would be interesting to determine the maximal $E_n$ suboperads
of $RU\scK^{(n)}$.

\bigskip
\noindent
{\bf Examples.} Theorem 2 allows us to construct many new examples of $E_n$
operads.
\begin{enumerate}
\item If we take $X=S^{n-1}$ with the antipodal $\Bbb Z/2$ action, then
$R_{\scP}X$ can be identified with an open suboperad of the topological
operad $|RU\scK^{(n)}|$ containing the closed suboperad $|\scP|$ as a
strong deformation retract. This follows from the fact that
$|\scK^{(n)}(2)|\cong S^{n-1}$.
\item Consider the preoperad $F(\Bbb R^n,\mbox{-})$ of configuration spaces
in $\Bbb R^n$. Then $UF(\Bbb R^n,\mbox{-})=F(\Bbb R^n,2)$ has a hemispherical
cellular decomposition (see Berger \cite{Be}), and $R_\scP F(\Bbb R^n,2)$ is an
$E_n$ operad containing $F(\Bbb R^n,\mbox{-})$ as a subpreoperad.
\item Consider the little $n$-cubes operad $\scC_n$. Again by Berger \cite{Be},
$U\scC_n=\scC_n(2)$ has a hemispherical cellular decomposition. Then
$R_\scP\scC_n$ is an $E_n$ operad containing $\scC_n$ as a subpreoperad,
but {\bf not} as a suboperad (cf. example above).
\item Let $X$, $Y$ be spaces with free $\Bbb Z/2$ actions, having the homotopy
type of $S^{m-1}$, $S^{n-1}$ respectively, and with hemspherical cellular
decompositions. Let $\scP_m$, $\scP_n$, $\scP_{m+n}$ denote either the
triple $\scM_m$, $\scM_n$, $\scM_{m+n}$ or $\scK^{(m)}$, $\scK^{(n)}$,
$\scK^{(m+n)}$. Then the join $X*Y$ has the homotopy type of $S^{m+n-1}$,
with a free $\Bbb Z/2$ action and with a hemispherical cellular decomposition.
Hence $R_{\scP_{m+n}}(X*Y)$ is an $E_{m+n}$ operad containing the $E_m$ operad
$R_{\scP_m}X$ and the $E_n$ operad $R_{\scP_n}Y$ as suboperads. Moreover
$R_{\scP_{m+n}}(X*Y)$ is a kind of ``homotopy tensor product'' of the operads
$R_{\scP_m}X$ and $R_{\scP_n}Y$: if $R_{\scP_{m+n}}(X*Y)$ acts on a space
$Z$ then any element of $R_{\scP_m}X(k)$ determines a map $Z^k\to Z$, which is
a homotopy homomorphism with respect to the $R_{\scP_n}Y$ actions on $Z^k$ and
$Z$, and similarly with the roles of $R_{\scP_m}X$ and $R_{\scP_n}Y$ reversed
(cf. \cite{BV} and \cite{BV2} and the discussion of tensor products of operads
below.)
\item Let $Y_1$, $Y_2$, \dots, $Y_n$ be arbitrary contractible spaces. Then
by iteration of the previous example, the iterated join
$$X=(Y_1\times\Bbb Z/2)*(Y_2\times\Bbb Z/2)*\dots*(Y_n\times\Bbb Z/2)$$
has the homotopy type of $S^{n-1}$ with free $\Bbb Z/2$ action and a
hemispherical cellular decomposition. Then $R_{\scP_n}X$ is an $E_n$ operad
containing the $E_1$ operads $R_{\scP_1}Y_i$ as suboperads, and $R_{\scP_n}X$
can be regarded as a homotopy tensor product of these suboperads.
\end{enumerate}

\bigskip
I now return to the proof of part (iii) of Theorem 1.
As I mentioned previously, this is accomplished via a more structured
right adjoint construction.

\bigskip
\noindent
{\bf Definition.} An {\it $\ell$-truncated operad\/} is a collection of
topological spaces $\{\scA(i)\}_{0\le i\le\ell}$ (with $\scA(0)=*$) with a
$\Sigma_i$ action on $\scA(i)$ and with compositions
$$\scA(k)\times\scA(i_1)\times\scA(i_2)\times\dots\scA(i_k)
\longrightarrow\scA(i_1+i_2+\dots+i_k)$$
defined whenever the spaces on both sides are defined, with these structures
satisfying all the relations required of an ordinary operad. The forgetful
functor
$$\mbox{Operads}\longrightarrow\mbox{$\ell$-truncated operads}$$
will be denoted $T_\ell$.

\bigskip
\noindent
{\bf Example.} A 1-truncated operad is just a monoid $\scA(1)$ (since
$\scA(0)=*$).

\bigskip
\noindent
{\bf Example.} A 2-truncated operad consists of the following data:
\begin{enumerate}
\item A monoid $\scA(1)$.
\item A $\Bbb Z/2$-space $\scA(2)$.
\item A left action
$$\scA(1)\times\scA(2)\longrightarrow\scA(2)$$
equivariant with respect to the $\Bbb Z/2$ action on both sides.
\item A right action
$$\scA(2)\times\left(\scA(1)\times\scA(1)\right)\longrightarrow\scA(2)$$
commuting with the left $\scA(1)$ action and equivariant with respect
to the $\Bbb Z/2$ action on both sides. ($\Bbb Z/2$ acts
on the left hand side via its action on $\scA(2)$ and by permuting
the factors of $\scA(1)\times\scA(1)$.)
\item A pair of maps
$$(d_1,d_2): \scA(2)\longrightarrow\scA(1)\times\scA(1)$$
which are equivariant with respect to the left action of $\scA(1)$
(diagonally on the right hand side), the right action of
$\scA(1)\times\scA(1)$, and the $\Bbb Z/2$ action.
\end{enumerate}

\bigskip
\noindent
{\bf Theorem 3.} (i) The forgetful functor
$$T_\ell:\mbox{Operads}\longrightarrow\mbox{$\ell$-truncated operads}$$
has a right adjoint
$$R_\ell:\mbox{$\ell$-truncated operads}\longrightarrow\mbox{Operads}$$
(ii) If $\scB$ is an operad with $\scB(1)=\{1_\scB\}$, then $R_2T_2\scB=RU\scB$
(where $U$ and $R$ are as in Theorem 1).

\bigskip
{\it Proof Sketch.\/} We describe the space $R_\ell\scA(k)$ as the space
of all orientations and labellings of the $(\ell-1)$-skeleton of the $(k-1)$-simplex
with vertices $\{1,2,\dots,k\}$. By an orientation of an $(i-1)$-subsimplex
($i\le\ell$) we mean a choice of total order of its vertices. By a
labelling of such a subsimplex we mean assigning a point in $\scA(i)$
corresponding to such a simplex. We impose an equivalence relation on
the orientations and labellings by specifying that changing the total
order of the vertices of an $i$-simplex by means of a permutation
in $\Sigma_i$ is equivalent to acting on its label by the same permutation.
We further impose the consistency condition on the orientations and
labellings of the subsimplices. The orientation and labelling on a face of
a subsimplex must be equivalent to the one obtained by restricting the
total order to the vertices and acting via the appropriate degeneracy.
(The total order on the vertices of the $(i-1)$-subsimplex determines a bijection
with the set $\{1,2,\dots,i\}$ and thus determines the appropriate degeneracy.)
More formally we can specify the points of $R_\ell\scA(k)$ as functions
defined on the set of pairs $(S,\lambda)$, where $S$ is a nonvoid subset of
$\{1,2,\dots,k\}$ of cardinality $\le\ell$ and $\lambda$ is a total order
on $S$, taking values in the disjoint union $\scA(1)\amalg\scA(2)\amalg\dots
\amalg\scA(\ell)$, and satisfying the conditions described above.

The action of $\Sigma_k$ on $R_\ell\scA(k)$ is via its action on the
vertices and orientations. The operad composition
$$R_\ell\scA(k)\times R_\ell\scA(i_1)\times R_\ell\scA(i_2)\times\dots
\times R_\ell\scA(i_k)\longrightarrow
R_\ell\scA(i_1+i_2+\dots+i_k)$$
can be described as follows. Consider an $(i-1)$ simplex with vertices
in $\{1,2,\dots,i_1+i_2+\dots+i_k\}$. Group them into blocks corresponding
to the factors on the left hand side. Then order the vertices within the
$j$-th block using the order in $R_\ell\scA(i_j)$. Order vertices in
different blocks using the order specified in $R_\ell\scA(k)$. If there
are $p$ different blocks of vertices and the blocks contain $q_1$, $q_2$,
\dots, $q_p$ vertices respectively, pick the labels specified in the product
space on the left hand side and multiply them together using the $\ell$-truncated operad composition
$$\scA(p)\times\scA(q_1)\times\scA(q_2)\times\dots\times\scA(q_p)
\longrightarrow\scA(q_1+q_2+\dots+q_p)$$

The unit of the adjunction $\scB\to R_\ell T_\ell\scB$ is given by appropriate
degeneracy maps. Part (ii) is obvious (which also proves Theorem 1(iii)).

\bigskip
\noindent
{\bf Remark.} The construction $R_1$ is due to Igusa \cite{I}. He calls this
construction the ``atomic operad generated by a monoid''.  We have
$$R_1M(k)=M^k$$
(where $M=\scA(1)$ is the given monoid, a.k.a. 1-truncated operad) with the
operad composition
$$M^k\times M^{i_1}\times M^{i_2}\times\dots\times M^{i_k}
\longrightarrow M^{i_1+i_2+\dots+i_k}$$
specified by the formula
\begin{eqnarray*}
\lefteqn{
\left(m_1,m_2,\dots,m_k;x_{11},\dots,x_{1i_1},x_{21},\dots,x_{2i_2},\dots,
x_{k1},\dots, x_{ki_k}\right)}\\
&\quad\mapsto
&\left(m_1x_{11},\dots,m_1x_{1i_1},m_2x_{21},\dots,m_2x_{2i_2},\dots,
m_kx_{k1},\dots, m_kx_{ki_k}\right)
\end{eqnarray*}
(It would be more accurate to call $R_1M$  the operad ``cogenerated by $M$'',
since this is a right adjoint construction and the forgetful functors $T_\ell$
also have left adjoints.)

\bigskip
\noindent
{\bf Remark.} The little $n$-cubes operads $\scC_n$ are 2-cogenerated, ie.
$$\scC_n\stackrel{\cong}{\longrightarrow}R_2T_2\scC_n.$$
To see this note that a $k$-fold configuration of little $n$-cubes in
$\scC_n(k)$ is determined by specifying each little $n$-cube in the
configuration, ie. by the vertex labels in $R_2T_2\scC_n(k)$. The labelling
of the edges in $R_2T_2\scC_n(k)$ specifies that each pair of little
$n$-cubes in the configuration specified by the vertex labels determines an
element of $\scC_n(2)$. This is equivalent to specifying that the interiors
of any two little $n$-cubes in the configuration specified by the vertex
labels are disjoint. Thus the configuration specified by the vertex labels
in $R_2T_2\scC_n(k)$ determines a unique element of $\scC_n(k)$. A similar
argument show that the linear isometries $E_\infty$ operad $\scL$ (cf.
\cite{BV}) is 2-cogenerated. Recall that the elements of $\scL(k)$ are linear
isometries $\left(\Bbb R^\infty\right)^k\to\Bbb R^\infty$, or equivalently a
$k$-tuple of linear isometries $\Bbb R^\infty\to\Bbb R^\infty$ whose images are
mutually orthogonal. The vertex labels of an element of $R_2T_2\scL(k)$
specify a $k$-tuple of linear isometries $\Bbb R^\infty\to\Bbb R^\infty$,
whereas the edge labels specify that their images are mutually orthogonal.
It follows that we have a homeomorphism of operads
$$\scL\stackrel{\cong}{\longrightarrow}R_2T_2\scL.$$

\bigskip
\noindent
I will now turn to an additional method of constructing $E_n$ operads, via
generalized tensor products.

\bigskip
\noindent
{\bf Definition.} Let $\scA$ and $\scB$ be operads acting on a topological space
$X$. We say that the actions of $\scA$ and $\scB$ are {\it interchangeable\/} if
for any $k$ and any element $\alpha\in\scB(k)$ the induced map 
$\alpha: X^k\longrightarrow X$ is a homomorphism of $\scB$-spaces, where $\scB$
acts on $X^k$ coordinatewise via its action on $X$. Equivalently for any 
$\beta\in\scB(\ell)$ we have
$$\alpha\cdot(\beta\times\beta\times\dots\times\beta)
(x_{\!\!\!\!\begin{array}{l}\scriptscriptstyle 1\le i\le k\\ [-5pt]\scriptscriptstyle 1\le j\le \ell\end{array}})
=\beta\cdot(\alpha\times\alpha\times\dots\times\alpha)\tau
(x_{\!\!\!\!\begin{array}{l}\scriptscriptstyle 1\le i\le k\\ [-5pt]\scriptscriptstyle 1\le j\le \ell\end{array}}),
$$
where $\tau:X^{k\ell}\to X^{k\ell}$ is the permutation which reorders the
coordinates of $X^{k\ell}$ from lexicographic to reverse lexicographic order.

\bigskip
\noindent
{\bf Definition.} Let $\scA\to\scC$ and $\scB\to\scC$ be operad maps. We say
that these operad maps are {\it interchangeable\/} if the following diagrams
commute for all $k$ and $\ell$:
$$
\xymatrix{\\ (*)\\ }\qquad\qquad
\xymatrix{
\scA(k)\times\scB(\ell)\ar[rr]^{id\times\Delta}\ar[d]^{\Delta\times id}
&&\scA(k)\times\scB(\ell)^k\ar[rr]
&&\scC(k)\times\scC(\ell)^k\ar[d]^{\mu}\\
\scA(k)^{\ell}\times\scB(\ell)\ar[d]^{\cong}
&& &&\scC(k\ell)\\
\scB(\ell)\times\scA(k)^{\ell}\ar[rr]
&&\scC(\ell)\times\scC(k)^{\ell}\ar[rr]_{\mu}
&&\scC(k\ell)\ar[u]^{\tau}
}
$$
Here $\mu$ denotes composition in the operad $\scC$ and $\tau$ denotes the
same permutation as in the preceding definition.

It is clear that if $\scC$ acts on a space $X$, then the induced actions by
$\scA$ and $\scB$ on $X$ are interchangeable. It was shown in \cite{BV} (cf. also
\cite{BV2}) that for any two operads $\scA$ and $\scB$ there is a pair
of operad maps $\scA\to\scA\otimes\scB$, $\scB\to\scA\otimes\scB$ which are
universal for interchangeable pairs of operad maps, ie. any other interchangeable
pair $\scA\to\scC$, $\scB\to\scC$ factors through a unique operad map
$\scA\otimes\scB\to\scC$. The operad $\scA\otimes\scB$ is called the {\it tensor
product\/} of $\scA$ and $\scB$.

In general analyzing the homotopy type of the tensor product of operads is
an intractable problem. However Dunn \cite{D} showed that the $n$-fold iterated
tensor product of little 1-cubes operads is an $E_n$ operad, more precisely
the suboperad of decomposable little $n$-cubes in $\scC_n$. We shall see below
another example of a tensor product of operads which can be shown to be $E_n$.
It is much easier to construct examples of $E_n$ operads which contain pairs
of interchangeable suboperads. We shall refer to such examples as generalized
tensor products.  More precisely if $\scA$ is an $E_m$ operad, $\scB$ is an
$E_n$ operad, and $\scC$ is an $E_{m+n}$ operad containing $\scA$ and $\scB$ as
interchangeable suboperads, we call $\scC$ a {\it generalized tensor product\/}
of $\scA$ and $\scB$.

\bigskip
\noindent
{\bf Definition.} A {\it partial acyclic orientation\/} of the complete graph
on the set of vertices $\{1,2,3,\dots,k\}$ is an assignment of direction to some
of the edges of the graph such that no directed cycles occur.  A 
{\it partial coloring\/} of the complete graph on $k$ vertices is an assignment of
colors to some of the edges of the graph from the countable set of colors
$\{1,2,3,\dots\}$.  The poset $\widehat{\scK}(k)$ has as elements pairs
$(\mu,\sigma)$, where $\mu$ is a partial coloring and $\sigma$ is a partial acyclic
orientation of the complete graph on $k$ vertices, with the condition that uncolored
edges are also unoriented and vice-versa. The order relation on
$\widehat{\scK}(k)$ is determined as follows: we say that
$(\mu_1,\sigma_1)\le (\mu_2,\sigma_2)$ if every uncolored unoriented edge
in $(\mu_1,\sigma_1)$ is also uncolored unoriented in $(\mu_2,\sigma_2)$, and for
any colored oriented edge $a\stackrel{i}{\longrightarrow}b$ in $(\mu_1,\sigma_1)$
the corresponding edge in $(\mu_2,\sigma_2)$ is either uncolored unoriented or has
either orientation and coloring $a\stackrel{j}{\longrightarrow}b$ with $j\ge i$
or $b\stackrel{j}{\longrightarrow}a$ with $j>i$.  The $n$-th filtration
$\widehat{\scK}^{(n)}(k)$ is the subposet of ${\scK}(k)$ where the colorings
are restricted to take values in the subset $\{1,2,3,\dots,n\}$.

The action of the symmetric group $\Sigma_k$ on $\widehat{\scK}(k)$ is via permutation
of the vertices. The composition
$$\widehat{\scK}(k)\times\widehat{\scK}(m_1)\times\widehat{\scK}(m_2)\times\dots\times\widehat{\scK}(m_k)
\longrightarrow\widehat{\scK}_(m_1+m_2+\dots+m_k)$$
assigns to a tuple of partial orientations and colorings in
$\widehat{\scK}(k)\times\widehat{\scK}(m_1)\times\widehat{\scK}(m_2)\times\dots\times\widehat{\scK}(m_k)$
the partial orientation and coloring obtained by subdividing the set of
$m_1+m_2+\dots+m_k$ vertices into $k$ blocks containing $m_1$, $m_2$, \dots, $m_k$
vertices respectively. The edges connecting vertices within the $i$-th block
are oriented and colored (or unoriented and uncolored) according to the given element
in $\widehat{\scK}(m_i)$.  The edges connecting vertices between blocks $i$ and $j$
are all oriented and colored (or not) according to the corresponding edge in the
given element of $\widehat{\scK}(k)$. It is easy to check that this specifies the
structure of a filtered operad on $\widehat{\scK}$ containing the complete graphs
operad $\scK$ as a filtered suboperad. Note however that $\widehat{\scK}^{(n)}$ is
not $E_n$: $\widehat{\scK}^{(n)}(k)$ is equivariantly contractible to the $\Sigma_k$
fixed point specified by the complete graph on $\{1,2,3,\dots,k\}$ with all its
edges unoriented and uncolored. We shall refer to $\widehat{\scK}$ and its filtrations
as the {\it augmented complete graphs operad\/}.

\bigskip
\noindent
{\bf Definition.} An $E_n$ operad $\scA$ is said to have an {\it augmented cellular
decomposition\/} over the complete graphs operad if
\begin{enumerate}
\item $\scA$ has a cellular decomposition over $\scK^{(n)}$.
\item $R_1\scA(1)$ has a cellular decomposition over $\widehat{\scK}^{(n)}$, where
$R_1$ is the right adjoint construction from monoids to operads of Theorem 3.
\item The adjunction map $\scA\to R_1\scA(1)$ is an imbedding onto the union of cells
indexed by $\scK^{(n)}$.
\end{enumerate}

\bigskip
\noindent
The basic example of an $E_n$ operad with an augmented cellular decomposition over
the complete graphs operad is the little $n$-cubes operad $\scC_n$. We specify the
cellular decomposition of $\left(R_1\scC_n(1)\right)(k)$ over $\widehat{\scK}^{(n)}(k)$
as follows. Given a partial orientation and coloring of the complete graph
on $\{1,2,3,\dots,k\}$ we define the corresponding cell to be the subspace of
$\scC_n(1)^k$ consisting of $k$-tuples of subcubes of the unit $n$-cube satisfying
\begin{enumerate}
\item For every oriented and colored edge $a\stackrel{i}{\longrightarrow}b$, the
interior of the $a$-th subcube must be separated from the interior of the $b$-subcube
by a hyperplane perpendicular to the $j$-th coordinate axis for some $j\le i$. If
$j=i$, then the $a$-th subcube is required to lie on the negative side of the
hyperplane, and the $b$-th subcube on the positive side.
\item If the edge joining $a$ to $b$ is unoriented and uncolored, no condition is
imposed on the relative positions of the $a$-th and $b$-th subcubes of the $k$-tuple,
eg. their interiors are allowed to intersect, even coincide.
\end{enumerate}
It is not difficult to check that this is a cellular decomposition satisfying the requirements
of the definition, using the arguments of \cite{Be}.

\bigskip
\noindent
{\bf Remark.} It may be worthwhile to note that requirements for an operad to have
an augmented cellular decomposition over the complete graphs operad are completely
antithetical to the requirement $\scA(1)=\{1_\scA\}$ for the construction of Theorem 1
to have good properties. Moreover while there are many examples of finite simplicial
$E_n$ operads satisfying $\scA(1)=\{1_\scA\}$, this is precluded for operads with
an augmented cellular decomposition.  Indeed the mere condition that 
$$(d_1,d_2):\scA(2)\longrightarrow \left(R_1\scA(1)\right)=\scA(1)^2$$
be injective precludes $\scA(k)$ from being a finite simplicial complex for any $k>0$.

To see this, pick any vertex $c\in\scA(2)$ and let $a=d_1(c)\in\scA(1)$. Then, assuming
that $\scA(1)$ is a finite complex, the  positive powers $a^m$ (with respect to the
monoid structure on $\scA(1)$) can't be all distinct. Suppose $a^m=a^{m+r}$. Consider
the element $c'=a^m c\, (a^{r-1}d_2(c),1)$, where we use the left action of $\scA(1)$
on $\scA(2)$ together with the right action of $\scA(1)\times\scA(1)$. Then
$$(d_1,d_2)(c') = \left(a^{m+k}d_2(c),a^m d_2(c)\right) = 
\left(a^m d_2(c), a^{m+k}d_2(c)\right) = (d_1,d_2)(c'\,\tau),$$
where $\tau$ is the transposition in $\Sigma_2$. It now follows that
$c'=c'\,\tau$, contradicting freeness of the $\Sigma_2$ action. Thus $\scA(1)$ can't
be a finite complex. It can be shown that for any $k\ge 2$, $\scA(k)$ can't be
a finite complex, by assuming the contrary, taking a vertex $c\in\scA(k)$, letting
$a=d_1d_1\dots d_1(c)\in\scA(1)$, noting that $a^m c$ can't all be distinct, and
arguing as above.

\bigskip
\noindent
{\bf Theorem 4.} (i) Let $\scA$ be an $E_m$ operad with an augmented cellular decomposition
over the complete graphs operad. Let $X$ be a space homotopy equivalent to $S^{n-1}$, with
a free $\Bbb Z/2$ action and a hemispherical cellular decomposition. Then there is a
generalized tensor product $\scA\gtensor R_{\scK^{(n)}}X$, an $E_{m+n}$ operad containing
$\scA$ and $R_{\scK^{(n)}}X$ as interchangeable suboperads.\newline
(ii) If $X=S^0$, then
$$\scA\gtensor R_{\scK^{(n)}}S^0 = \scA\otimes R_{\scK^{(n)}}S^0 = \scA\otimes\scM,$$
the actual tensor product of $\scA$ with $\scM$, the $A_\infty$ operad which acts on
(strict) monoids.

\bigskip
{\it Proof Sketch.\/} First consider the product operad $R_1\scA(1)\times R_{\scK^{(n)}}X$.
The elements of the $k$-th space of this operad can be described as the complete graph
on the set $\{1,2,\dots,k\}$, with the vertices being labelled by elements of $\scA(1)$
and the edges labelled by elements of $X$.

Now let us impose the quotient relation that if the endpoints of an edge are labelled
by a pair of elements in $\scA(1)\times\scA(1)$ which is in the image of $\scA(2)$,
then we are allowed to replace the $X$ label on that edge by any other element of $X$.
It is straightforward to check that the composition operation in $R_1\scA(1)\times R_{\scK^{(n)}}X$
passes to the resulting collection of quotient spaces.  Thus we obtain an operad we denote
$\scA\gtensor R_{\scK^{(n)}}X$.

There is an evident inclusion of operads $R_{\scK^{(n)}}X\subset\scA\gtensor R_{\scK^{(n)}}X$,
given by labelling all vertices in the complete graph by $1_{\scA}$.  There is also an inclusion
of operads $\scA\subset \scA\gtensor R_{\scK^{(n)}}X$ given by labelling all edges of the complete
graph by arbitrary elements of $X$ (which is well-defined by the quotient relation we have imposed).

These suboperads are easily seen to be interchangeable: if we take $\scB= R_{\scK^{(n)}}X$ in
diagram (*) and chase an element $(\alpha, \beta)\in\scA(k)\times \left(R_{\scK^{(n)}}X\right)(\ell)$,
around the two sides of the diagram, we obtain the same labelling of the complete graph on $k\ell$ vertices
either way.  For if we identify $\{1,2,\dots,k\ell\}$ with the product 
$\{1,2,\dots,k\}\times\{1,2,\dots,\ell\}$ via lexicographic ordering on the latter, then the labelling
of the complete graph obtained by going along the top and right of the diagram assigns to vertex $(i,j)$
the label $\alpha_i=\alpha(*\times *\times\dots 1_{\scA}\times\dots\times *)$, where $1_{\scA}$
is in the $i$-th place. The label assigned to the edge joining vertices $(i_1,j_1)$ and $(i_2,j_2)$ is
$\beta_{j_1j_2}$ provided that $i_1=i_2$ and $j_1\ne j_2$, otherwise it is unspecified (which is
allowed by the quotient relation imposed above). If we go around the other side of the diagram, we obtain
the same labelling on the vertices. The labelling assigned to the edge joining vertices $(i_1,j_1)$ and
$(i_2,j_2)$ is $\beta_{j_1j_2}$ provided that $j_1\ne j_2$. However the quotient relation allows
us to discard those labels when $i_1\ne i_2$, so we obtain the same labelling as before.

We show that $\scA\gtensor R_{\scK^{(n)}}X$ is an $E_{m+n}$ operad by displaying a cellular decomposition
over $\scK^{(m+n)}$. Given an element $\lambda=(\mu,sigma)\in\scK^{(m+n)}(k)$, the corresponding cell in
$\left(\scA\gtensor R_{\scK^{(n)}}X\right)(k)$ can be described as follows.  Consider the cell in $\scA(1)^k$
specified by $(\mu_1, \sigma_1)\in\widehat{\scK}^{(m)}(k)$ specified by deleting directions and colorings of those
edges of the complete graph whose colors are $>m$.  This specifies the allowable labellings of vertices of the
complete graph of the elements in $\left(\scA\gtensor R_{\scK^{(n)}}X\right)(k)$.  Now for each edge whose
color is $j>m$, let the allowable labellings of that edge be elements of $X$ specified by the given orientation of
that edge and color $j-m$.  Then this cell in $\left(\scA\gtensor R_{\scK^{(n)}}X\right)(k)$ is homeomorphic
to the pushout of the following diagram
$$
\xymatrix{
\bigcup_{\lambda'_1<\lambda_1}C_{\lambda'_1}\times D_{\lambda-\lambda_1}\ar@{^{(}->}[r]\ar[d]
&C_{\lambda_1}\times D_{\lambda-\lambda_1}\\
\bigcup_{\lambda'_1<\lambda_1}C_{\lambda'_1}\times D_{\lambda-\lambda'_1}
}
$$
where $\lambda_1(\mu_1, \sigma_1)\in\widehat{\scK}^{(m)}(k)$, and $C_{\lambda_1}$ is the corresponding
cell in $\scA(1)^k$.  $D_{\lambda-\lambda_1}\subset X^r$, where $r$ denotes the number of edges which are not
assigned a coloring and orientation by $\lambda_1$, and which therefore specify elements of $X$ as indicated above.
Similarly for $\lambda'_1<\lambda_1\in\widehat{\scK}^{(m)}(k)$.  This pushout is contractible since 
$C_{\lambda_1}\times D_{\lambda-\lambda_1}$ is contractible, the horizontal arrow is a cofibration, and
the vertical arrow is an equivalence.

Part (ii) follows by comparing with the explicit description of $\scA\otimes\scM$ given in \cite{BV2}.

\bigskip
\noindent
{\bf Remark.} In general $\scA\otimes R_{\scK^{(n)}}X$ is obtained from $R_1\scA(1)\times R_{\scK^{(n)}}X$
by imposing a finer equivalence than that used to define $\scA\gtensor R_{\scK^{(n)}}X$.  Although one could define
``cells'' in the resulting operad similarly to the above, it is difficult to check under what circumstances they
are contractible.  One would certainly need to impose further conditions on $X$.

\bigskip
\noindent
{\bf Theorem 5.} Let $\scA$ be an $E_m$ operad and let $\scB$ be an $E_n$ operad, both
with augmented cellular decompositions over the complete graphs operad. Then there is a
generalized tensor product $\scA\gtensor \scB$, an $E_{m+n}$ operad containing
$\scA$ and $\scB$ as interchangeable suboperads.

\bigskip
{\it Proof Sketch.\/} First consider the product operad $R_1\scA(1)\times R_1\scB$.
The elements of the $k$-th space of this operad can be described as the complete graph
on the set $\{1,2,\dots,k\}$, with the vertices being labelled by elements of $\scA(1)\times\scB(1)$.
We define $\scA\gtensor \scB$ to be the suboperad consisting of those elements such that for
any pair of vertex labels $(\alpha_1,\beta_1)$, $(\alpha_2,\beta_2)$ at least one of $(\alpha_1,\alpha_2)$
or $(\beta_1,\beta_2)$ is in $\scA(2)\subset\scA(1)\times\scA(1)$ or $\scB(2)\subset\scB(1)\times\scB(1)$.

Then $\scA\subset R_1\scA\times R_1\{1_{\scB}\}$ and $\scB\subset R_1\{1_{\scA}\}\times R_1\scB$
are both suboperads of $\scA\gtensor \scB$, and are easily seen to be interchangeable.

Now for each $(\mu,\sigma)\in\scK^{(m+n)}(k)$ define $(\mu_1,\sigma_1)\in\widehat{\scK}^{(m)}(k)$
to be the partial coloring and orientation obtained from $(\mu,\sigma)$ by removing the coloring and
orientation of any edge colored by a color $>m$. Similarly define $(\mu_2,\sigma_2)\in\widehat{\scK}^{(n)}(k)$
to be the partial coloring and orientation obtained from $(\mu,\sigma)$ by removing the coloring and
orientation of any edge colored by a color $\le m$ and changing the color of all the remaining edges according
to the rule $\i\mapsto i-m$.  Now let the cell corresponding to $(\mu,\sigma)$ be
$C_{(\mu_1,\sigma_1)}\times D_{(\mu_2,\sigma_2)}\subset \scA(1)^k\times\scB(1)^k$.
It is easy to check that this specifes a cellular decomposition of $\scA\gtensor \scB$ over
$\scK^{(m+n)}$.

\bigskip
\noindent
{\bf Remarks.} (i) If in the above construction we take $\scA$ and $\scB$ to be the little $m$-cubes,
resp. little $n$-cubes operads, then $\scA\gtensor \scB$ is the little $(m+n)$-cubes operad.\newline
(ii) It is easy to see that $\scA\gtensor \scB$ actually is an $E_{m+n}$ operad 
with an augmented cellular decomposition over the complete graphs operad.  Thus the construction can
be iterated.

\end{document}